\documentclass[3p,preprint,11pt]{elsarticle}

\usepackage[utf8]{inputenc}
\usepackage{amsmath}
\usepackage[english]{babel}
\usepackage{amssymb}
\usepackage{amsfonts}
\usepackage{bm}
\usepackage{graphicx}
\usepackage{amsmath}
\usepackage{amsthm}
\usepackage{xcolor}
\usepackage{parskip}
\usepackage{multicol}
\usepackage{float}
\usepackage{rotating}
\usepackage{rotating}

\usepackage{chngcntr}
\usepackage{apptools}

\usepackage{amssymb}

\journal{}

\usepackage{graphicx,enumitem}

\usepackage{subfigure}

 \theoremstyle{definition}
\newtheorem{theorem}{Theorem}[section]
\newtheorem{cor}{Corollary}[section]
\newtheorem{lemma}{Lemma}[section]
\newtheorem{proposition}{Proposition}[section]
\newtheorem{definition}{Definition}[section]

\newtheorem{ejemplo}{Example}[section]

\usepackage{natbib}

\numberwithin{equation}{section}
\numberwithin{figure}{section}
\numberwithin{table}{section}

\usepackage[colorlinks=true,linkcolor=red, citecolor=blue]{hyperref}

\def\R{\mathbb{R}}
\def\tr{\text{trace}}

\begin{document}

\begin{frontmatter}

\title{The dimple problem related to space-time modeling under the Lagrangian framework}

\author[1]{Alfredo Alegr\'ia\footnote{Corresponding author. Email: alfredo.alegria.jimenez@gmail.com}}
\author[2]{Emilio Porcu}
\address[1]{Departamento de Matem\'atica, Universidad T\'ecnica Federico Santa Mar\'ia, Valpara\'iso, Chile}
\address[2]{School of Mathematics and Statistics,   Newcastle University, Newcastle,  United Kingdom}

\begin{abstract} 
 Space-time covariance modeling under the Lagrangian framework has been especially popular to study atmospheric phenomena in the presence of transport effects, such as prevailing winds or ocean currents, which are incompatible with the assumption of full symmetry. In this work, we assess the dimple problem (Kent et al., 2011) for covariance functions coming from  transport phenomena. We work under two important cases: the spatial domain can be either the $d$-dimensional Euclidean space $\mathbb{R}^d$ or the spherical shell of $\mathbb{R}^d$. The choice is relevant for the type of metric chosen to describe spatial dependence. In particular, in Euclidean spaces, we work under very general assumptions with the case of radial symmetry being deduced as a corollary of a more general result.  We illustrate through  examples that, under this framework,  the dimple  is  a natural  and physically interpretable property.  
\end{abstract}

\begin{keyword}
Covariance functions   \sep Random fields \sep Rotation group \sep Transport effects.
\end{keyword}

\end{frontmatter}

\section{Introduction}

Space-time geostatistics deals mainly with the second order properties of  random fields   evolving temporally over a given  spatial domain. In particular, covariance functions describe the interactions between spatial and temporal components, and they are  crucial for both estimation and prediction. A thorough description of the properties of space-time covariance functions is given by  \citet{GGG07}. Throughout this work, we consider  zero mean space-time  fields, $\{Z(\bm{x},t): \bm{x} \in {\cal D}, t \in \R\}$,  where $\bm{x}$ is a spatial location in the region $\cal{D}$ and $t$ is a temporal instant.  The covariance function associated to $Z(\bm{x},t)$ is defined through the mapping $K: ({\mathcal D} \times  \R)^2 \to \mathbb{R}$ given by $K(\bm{x},t,\bm{y},t') = {\rm cov}\left\{ Z(\bm{x},t), Z(\bm{y},t') \right \}$, for $\bm{x}, \bm{y} \in {\cal D}$ and $t, t' \in \R$. \citet{GGG07} describe space-time covariances when ${\cal D}= \mathbb{R}^d$, where $d$ denotes a positive integer. We instead refer to \citet{berg2016schoenberg} and \citet{doi:10.1080/01621459.2015.1072541}    for the case of  fields evolving temporally over ${\cal D}=\mathbb{S}^{d-1} = \{\bm{x}\in\mathbb{R}^d, \|\bm{x}\|=1\}$, with $\|\cdot\|$ denoting the Euclidean distance. In particular, the case  $\mathcal{D}=\mathbb{S}^2$ can be used in statistical practice to represent planet Earth.

The work by \citet{Kent01062011} has brought attention to the so-called dimple problem: for some classes of space-time covariance functions,  a dimple is present if $Z(\bm{x}_{\text{here}},t_{\text{now}})$ is more correlated with $Z(\bm{x}_{\text{there}},t_{\text{tomorrow}})$ than with  $Z(\bm{x}_{\text{there}},t_{\text{now}})$.  The authors establish conditions for the presence of dimple in the  Gneiting  covariance  \cite{10.2307/3085674} and argue that the dimple is a counterintuitive property for modeling space-time data since it \textit{contradicts a natural monotonicity requirement} of the covariance.  Additional works related to the dimple effect are \citet{cuevas},  \citet{heidok22193}, \citet{2015arXiv150501243H} and \citet{Mosammam2014}.

 In this paper, we study the dimple problem for space-time covariances coming from  transport effect (or Lagrangian) models, which are a popular alternative to analyze atmospheric phenomena in the presence of flow, such as prevailing winds or waves. For a detailed discussion on transport effect   models the reader is referred to  \cite{cressie2011statistics,GGG07,JGRD:JGRD780} and the extensive list of references therein. Recently, Fiedler \cite{heidok22193} has provided an interesting example of precipitation data for three German cities, where  the empirical covariance exhibits a dimple-like behavior influenced by prevailing westerly winds. In addition, \citet{heidok22193}  performs simulation experiments with transport effect models and shows that the dimple can be relevant in  statistical applications.  Indeed, under the Lagrangian framework, the dimple effect is an expected and physically interpretable property.  
 
  We devote separate expositions to the cases where ${\cal D}$ can be either the $d$-dimensional Euclidean space or the unit sphere $\mathbb{S}^{d-1}$ of $\mathbb{R}^d$. The choice of the spatial domain has a crucial effect on the metric describing the distances between any pair of points, hence on the structure of the covariance function and its mathematical representation as well. Furthermore, previous literature on the dimple effect in $\mathbb{R}^d$ is based on the assumption of stationarity and isotropy, namely that the covariance function is radially symmetric in the spatial argument. We work under more general assumptions and show that the dimple effect can be analyzed under  general frameworks.   On the other hand, we could not find in previous literature any mathematical formulation of a transport effect model when the space is the spherical shell of $\mathbb{R}^d$. The classical formulation described in \cite{GGG07} would not apply here, because the curvature of the sphere  must be taken into account. We have been able to propose such architecture through the use of spatial rotations, and this is the crux for constructing the associated covariance, and then inspecting the dimple problem.

 The article is organized as follows. Section 2 studies the dimple problem for  transport effect models when ${\cal D}=\mathbb{R}^d$. Our characterization is then discussed under additional assumptions on the   transport directions, as well as under the assumption that the covariance generating the Lagrangian model is spatially isotropic. Section 3 is devoted to  fields on spheres across time where we introduce a new transport effect model, respecting the spherical curvature, and characterize the related dimple property.   We provide several examples where the dimple naturally arises as a result of  the transport dynamic.  In Section 4  we give some conclusions.

\section{Dimple effect for transport phenomena on $\mathbb{R}^d \times \mathbb{R}$}

We  start by introducing the definition of dimple for stationary fields  $Z(\bm{x},t)$ on $\R^d \times \R$, for which the covariance function is represented by a continuous mapping $C:\R^d \times \R \to \R$, that is, ${\rm cov}\{Z(\bm{x},t),Z(\bm{y},t')\}= C(\boldsymbol{h}, u)$, where we use the  notation $\boldsymbol{h} \in \R^d$ for the spatial lag $\bm{x}-\bm{y}$, and $u \in \R$ for the temporal lag $t-t'$.   

\begin{definition}	
\label{def_dimple}
Let $Z(\bm{x},t)$ be a weakly  stationary random field on $\R^d \times \R$, with stationary covariance $C: \R^d \times \R \to \R$. Then, $C$ has a dimple along the temporal lag $u \in \R$ if there exists two sets, $\Lambda _1 \subset \mathbb{R}^d$ and $\Lambda_2  \subset  \R^d \setminus \Lambda_1$, with $\Lambda_1$ containing the origin, such that the following conditions hold:
\begin{itemize}
\item[(1)] For  fixed $\bm{h}_0 \in \Lambda_1$, the mapping $u \mapsto C(\bm{h}_0,u)$, restricted to   $u\geq 0$, has a local maximum at $u=0$. 
\item[(2)] For  fixed $\bm{h}_0 \in \Lambda_2$, the mapping $u \mapsto C(\bm{h}_0,u)$, restricted to   $u\geq 0$, has a local minimum at $u=0$.
\end{itemize}
\end{definition}
Here,  $\Lambda_1$ and  $\Lambda_2$  can be interpreted as the sets of spatial lags where dimple is absent and present, respectively.  For example, the dimple effect could occur in a determined set of spatial directions, but not in another one. 

Since we are not assuming  time symmetry, $C(\bm{h},u) \neq  C(\bm{h},-u)$ in general.  Thus, Definition \ref{def_dimple} can be considered as a \textit{right-dimple} effect, by noting that $u \mapsto C(\bm{h}_0,u)$ is only studied for  non negative temporal lags. A similar  \textit{left-dimple} definition can be introduced by studying the mapping $u \mapsto C(\bm{h}_0,u)$ in $u\leq 0$, with analogous conditions and  interpretation. However,  a straightforward calculus, based on the identity $C(\bm{h},u)=C(-\bm{h},-u)$, shows that a right-dimple occurs, with  sets $\Lambda_i$, for  $i=1,2$,  if, and only if, a left-dimple occurs, with sets $\widetilde{\Lambda}_i=\{\bm{x}\in\mathbb{R}^d, -\bm{x}\in\Lambda_i\}$, for $i=1,2$.  Therefore, it is not worth making the distinction between them and we simply call it a dimple effect.

In particular, if $C$ is radially symmetric in space and symmetric in time, that is, if it depends on $\bm{h}$ and $u$ only through $\|\bm{h}\|$ and $|u|$, we obtain an equivalent definition to that  introduced by \citet{Kent01062011}. For such a specific case,  $\Lambda_1$  can be taken as the ball $\Lambda_1=\{\bm{h}\in\mathbb{R}^d, \|\bm{h}\| \leq  L \}$, for some  $L\geq 0$, whereas that $\Lambda_2$ can be taken as its complement. Thus,  Definition \ref{def_dimple} can be adapted as follows.
\begin{definition}	
\label{def_dimple2}
Let $Z(\bm{x},t)$ be a weakly stationary random field on $\R^d \times \R$, with covariance  $C: \R^d \times \R \to \R$ being radially symmetric in space and symmetric in time.  Then, $C$ has a dimple along the temporal lag $u \in \R$ if there exists  $L\geq 0$ such that the following conditions hold:
\begin{itemize}
\item[(1)] For  fixed $\|\bm{h}_0\| \leq  L$, the mapping $u \mapsto C(\bm{h}_0,u)$ has a local maximum at $u=0$.
\item[(2)] For  fixed $\|\bm{h}_0\|> L$, the mapping $u \mapsto C(\bm{h}_0,u)$ has a local minimum at $u=0$.
\end{itemize}
\end{definition}
The original definition introduced by \citet{Kent01062011} involves conditions on the  increasing or decreasing behavior of $u \mapsto C(\bm{h}_0,u)$ on the positive real line. The connection with Definition \ref{def_dimple2} is clear, since condition (1) is related to a decreasing behavior of $u \mapsto C(\bm{h}_0,u)$ in a right neighborhood of $u=0$, whereas that condition (2) is associated to an increasing local behavior of this mapping.

Our goal is to study the presence of dimple  for transport effect models. For a stationary and merely spatial Gaussian field $Y(\bm{x})$ on $\R^d$, with covariance function $C_S:\R^d \to \R$, and a random vector $\bm{V}$ in $\R^d$, we define a space-time field $Z(\bm{x},t)$ with transport effect according to 
\begin{equation}
\label{field1}
 Z(\bm{x},t)= Y(\bm{x}- t \boldsymbol{V}), \qquad (\bm{x},t) \in \R^d \times \R.
 \end{equation}
Equation (\ref{field1}) represents a space-time dynamic in which an entire spatial field  moves time-forward in some random  direction determined by the  velocity vector $\bm{V}$. Straightforward calculus show that the covariance function $C$ associated to $Z(\bm{x},t)$ is stationary and has expression
\begin{equation}\label{transporte}
C(\bm{h},u) = {\rm E}\{ C_{S}(\bm{h}- u\bm{V}) \}, \qquad  (\bm{h},u)\in \mathbb{R}^d\times \mathbb{R},
\end{equation}
where   expectation $\rm{E}\{\cdot\}$ is taken with respect to the $d$-dimensional random vector $\bm{V}$, whose   characteristic function is denoted as $\varphi_V (\bm{\eta}) = {\rm E} \{  \exp (\bm{\imath} \bm{V}^\top \bm{\eta} )  \}$, where $\bm{\imath}=\sqrt{-1}\in\mathbb{C}$ and $\top$ is the transpose operator.   We refer equivalently to $Z(\bm{x},t)$ in (\ref{field1}) or to the related covariance $C$ in (\ref{transporte}) as trasport effect model or Lagrangian framework. \citet{GGG07} argue that the choice of the random velocity vector in (\ref{transporte}) 
should be justified on the basis of physical considerations.  In general,  expression (\ref{transporte})   generates models that are neither  radially symmetric in space nor temporally symmetric.

We now characterize the dimple property, according to Definition \ref{def_dimple}, for the transport covariance  (\ref{transporte}).        We first pay attention to the simplest scenario, where the velocity vector $\bm{V}$ is constant and nonzero. In this case, the model is referred to frozen field by \citet{JGRD:JGRD780} and  it is typically used to represent a prevailing wind along a given direction. The covariance (\ref{transporte}) reduces to $C(\bm{h},u) = C_{S}(\bm{h}- u\bm{V})$ and the presence of dimple  can be shown by direct inspection.  In fact, for any $\bm{h}\neq \bm{0}$ being parallel to $\bm{V}$, the  maximum value of  the mapping $u\rightarrow C(\bm{h},u)$ is not reached at origin, and the dimple obviously occurs.  On the other hand, for  spatial lags  $\bm{h}\neq \bm{0}$ that are orthogonal to $\bm{V}$, the dimple effect does not occur. 

\begin{figure}
\centering
\includegraphics[scale=0.25]{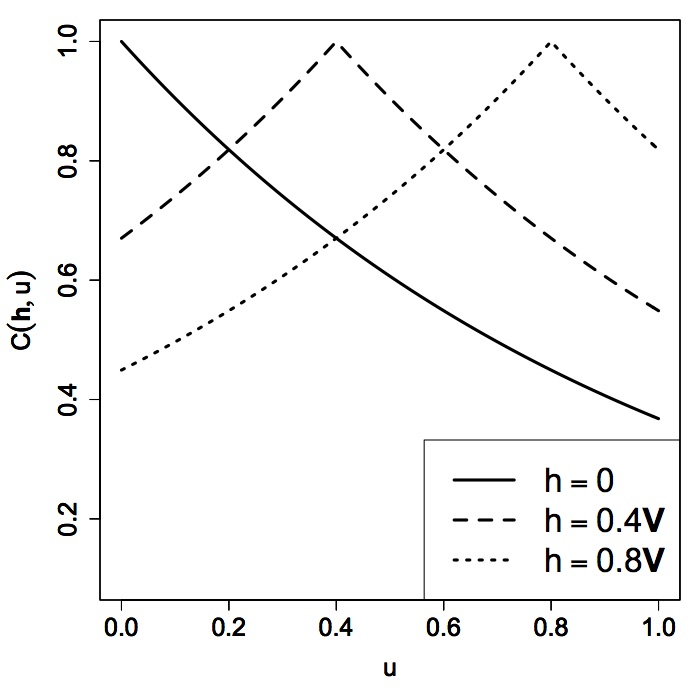}  \hspace{2cm}  \includegraphics[scale=0.25]{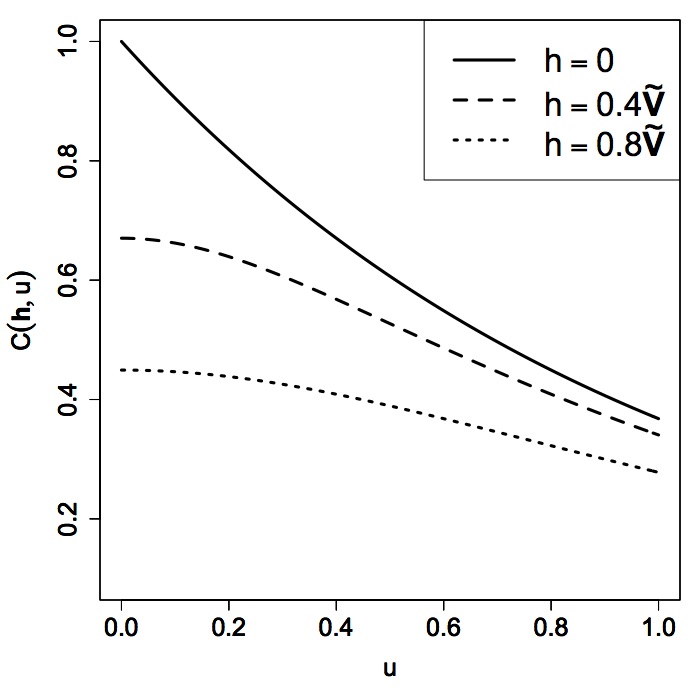}
\caption{The mapping $u \mapsto C(\bm{h},u)$ for the transport covariance $C(\bm{h},u) = C_{S}(\bm{h}- u\bm{V})$ on $\mathbb{R}^2\times \R$, with  deterministic direction $\bm{V}=(1,0)^\top$ and   exponential spatial covariance $C_S$. We fix spatial lags proportional to $\bm{V}$ (left), where dimple effect occurs, as well as spatial lags proportional to $\widetilde{\bm{V}}=(0,1)^\top$ (right), where dimple effect does not occur.}
\label{not_isotropic}
\label{constant}

\end{figure}

\begin{ejemplo}
Consider $d=2$, the  deterministic  direction $\bm{V}=(1,0)^\top$ and the exponential spatial covariance $C_S(\bm{h})=\exp(-\|\bm{h}\|)$.  Let  $\widetilde{\bm{V}}=(0,1)^\top$,  which is orthogonal to $\bm{V}$. Figure \ref{constant} illustrates the mapping $u\rightarrow C(\bm{h}_0,u)$,  where dimple effect is present for   spatial lags  that are proportional to $\bm{V}$, but it does not occur for spatial lags that are proportional to $\widetilde{\bm{V}}$.
\end{ejemplo}

\begin{figure}
\centering
\includegraphics[scale=0.25]{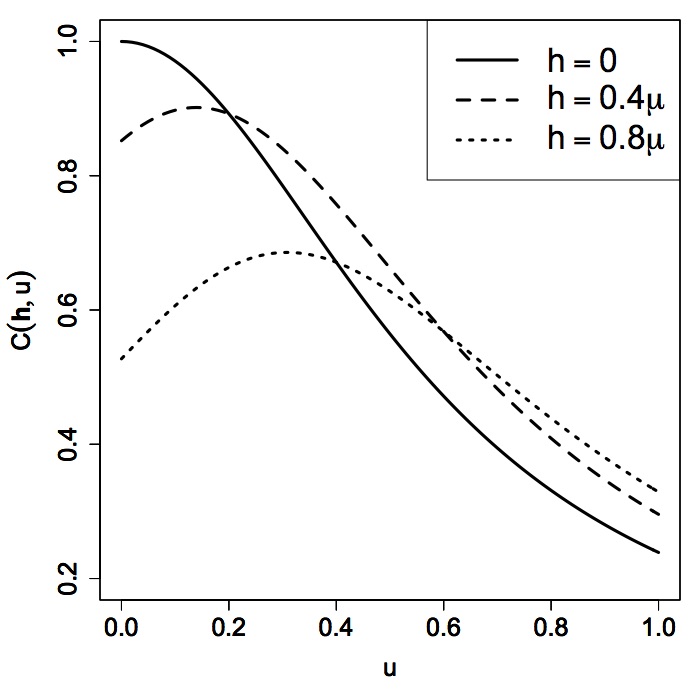}  \hspace{2cm}  \includegraphics[scale=0.25]{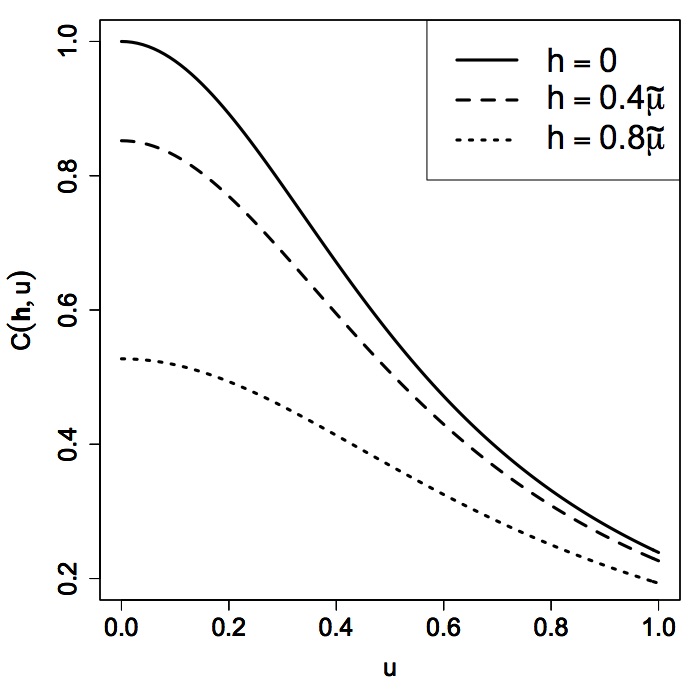}
\caption{The mapping $u \mapsto C(\bm{h},u)$ for the transport covariance $C(\bm{h},u) = {\rm E}\{C_{S}(\bm{h}- u\bm{V})\}$ on $\mathbb{R}^2\times \R$, with  $\bm{V}\sim \mathcal{N}_2(\bm{\mu},I_2)$, where $\bm{\mu}=(1,0)^\top$, and   Gaussian spatial covariance $C_S$. We fix spatial lags proportional to $\bm{\mu}$ (left), where dimple effect occurs, as well as spatial lags proportional to $\widetilde{\bm{\mu}}=(0,1)^\top$ (right), where dimple effect does not occur.}
\label{gauss}
\end{figure}

Another interesting   case arise when $\bm{V}\sim\mathcal{N}_d(\bm{\mu},\Sigma)$, where $\mathcal{N}_d$ denotes the $d$-dimensional Gaussian distribution.   \citet{schlather2010some} provides a closed form expression for the transport covariance under Gaussian distributed velocity vector.  The following example illustrates the dimple effect in this scenario.
\begin{ejemplo}
Consider $d=2$, $\bm{\mu}=(1,0)^\top$, $\Sigma=I_2$, with $I_d$ denoting the $d\times d$ identity matrix,  and the Gaussian spatial covariance $C_S(\bm{h})=\exp(-\|\bm{h}\|^2)$. The resulting transport covariance (see \citet{schlather2010some}) is given by 
$$C(\bm{h},u) = \frac{1}{1+2u^2}\exp\left( -\frac{(h_1-u)^2+h_2^2}{1+2u^2} \right), \qquad \bm{h}=(h_1,h_2)^\top\in\R^2, u\in\R.$$
Naturally, we expect the occurrence of dimple along  spatial lags that are proportional to $\bm{\mu}$, which is the main direction of the transport dynamic, and no dimple for   spatial lags that are proportional to  $\widetilde{\bm{\mu}}=(0,1)^\top$, which is orthogonal to $\bm{\mu}$  (see Figure \ref{gauss}).
\end{ejemplo}

The dimple effect can also occur  in circumstances where there is no prevailing directions. The presence of dimple is not directly obtained under this setting.   The rest of this section focuses on such cases. In particular, we consider   symmetrically distributed   vectors $\bm{V}$, that is, when $\bm{V}$ and $-\bm{V}$ have the same probability law.   This choice ensures that $C$ is temporally symmetric, but not necessarily radially symmetric in space. Dimple effects for the radial case will be deduced as corollary.    
 
 Before we state the main result of this section, we must introduce some notation. For a function $g:\mathbb{R}^d\rightarrow \mathbb{R}$, we denote $\nabla g(\bm{x} )= (\partial g(\bm{x})/ \partial x_i)_{i=1}^d$ the gradient vector associated to $g$. Also, we define $\nabla^2 g(\bm{x}) = \left[ \partial^2 g(\bm{x})/ \partial x_i \partial x_j \right]_{i,j=1}^d$ as its Hessian matrix. Accordingly, we define $\mathcal{H}_V=\nabla^2 \varphi_V (\bm{x})|_{\bm{x}=\boldsymbol{0}}$ as the symmetric  Hessian matrix of $\varphi_V$ evaluated  the origin.  The following result gives a characterization of the dimple effect, under symmetrically distributed directions, in terms of the second order derivatives of $C_S$ and $\varphi_V$.

\begin{theorem} 
\label{teorema1}
Let $C_S: \R^d \to \R$ be a continuous covariance function, being twice differentiable on
$\mathbb{R}^d\setminus \{\bm{0}\}$. Let $\bm{V}$ be symmetrically distributed on $\R^d$ 
 with characteristic function being twice differentiable  at origin. Let  $F:\mathbb{R}^d\setminus \{\bm{0}\}\rightarrow \mathbb{R}$ be defined as
\begin{equation} 
\label{aporte}
F(\bm{h}) = \tr\{ \mathcal{H}_V \ \nabla^2 C_S(\bm{h}) \}, \qquad \bm{h} \in \R^d.
\end{equation}
Then, the  covariance (\ref{transporte}) has a dimple if, and only if, there exists two sets  $\Lambda_1$  and $\Lambda_2$, according to Definition \ref{def_dimple},  such that   $F$ is positive in $\Lambda_1\setminus \{\bm{0}\}$,  and   negative  in $\Lambda_2$. 
\end{theorem}

Note that this result does not require differentiability conditions at origin for the spatial covariance $C_S$.   On the other hand, if $F$ is negative over each point of its domain,  the transport covariance $C$ has a dimple in the temporal lag with $\Lambda_1=\{\bm{0}\}$ and $\Lambda_2=\mathbb{R}^d\setminus\{\bm{0}\}$. In this case, the dimple arises immediately whenever the spatial separation $\bm{h}$ is nonzero.  The proof of Theorem \ref{teorema1} is deferred to  Appendix \ref{app1}.

\begin{ejemplo}
Consider $d=2$ and $\bm{V}$ a dichotomic random vector with distribution $\Pr(\bm{V} = \pm \bm{\xi}) = 0.5$, where $\bm{\xi} \in \mathbb{R}^2\setminus\{\bm{0}\}$ is a fixed vector. In this case, the transport model  is given by $C(\bm{h},u) = \{ C_S(\bm{h}-u\bm{\xi}) + C_S(\bm{h}+u\bm{\xi}) \}/2$, for $(\bm{h},u) \in \R^2 \times \R.$ Note that  $\varphi_V (\bm{\eta}) = \cos(   \bm{\eta}^\top\bm{\xi} )$ and  $\mathcal{H}_V = - \bm{\xi}  \bm{\xi}^\top$. Consider  a Cauchy spatial covariance $C_S(\bm{h}) = (1+\|\bm{h}\|^2)^{-1}$ and $\bm{\xi}=(1,1)^\top$. A direct calculation shows that the function $F(\bm{h})$ defined in Theorem \ref{teorema1} is given by $ F(\bm{h}) = 4(1- \|\bm{h}\|^2-4h_1h_2)  / (1+\|\bm{h}\|^2)^3$, where $\bm{h}=(h_1,h_2) \in \R^2$. We have that $F$ is negative along the points of the form  $\bm{h}=\kappa \bm{\xi}$, for $|\kappa| > 1/\sqrt{6}$,  and positive along  the directions that are orthogonal   to $\bm{\xi}$. Then, we should   expect the presence of dimple in the directions proportional to $\bm{\xi}$ (see Figure \ref{not_isotropic}).  Note that, unlike the previous examples and due that condition $|\kappa| > 1/\sqrt{6}$ is required, the dimple effect does not appear immediately for any $\bm{h}\neq \bm{0}$ proportional to $\bm{\xi}$.
\end{ejemplo}

\begin{figure}
\centering 
\includegraphics[scale=0.25]{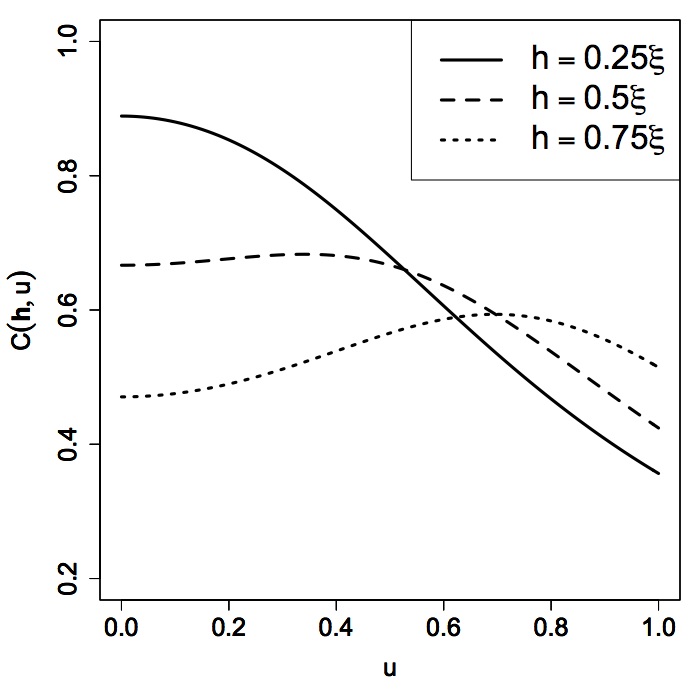}   \hspace{2cm} \includegraphics[scale=0.25]{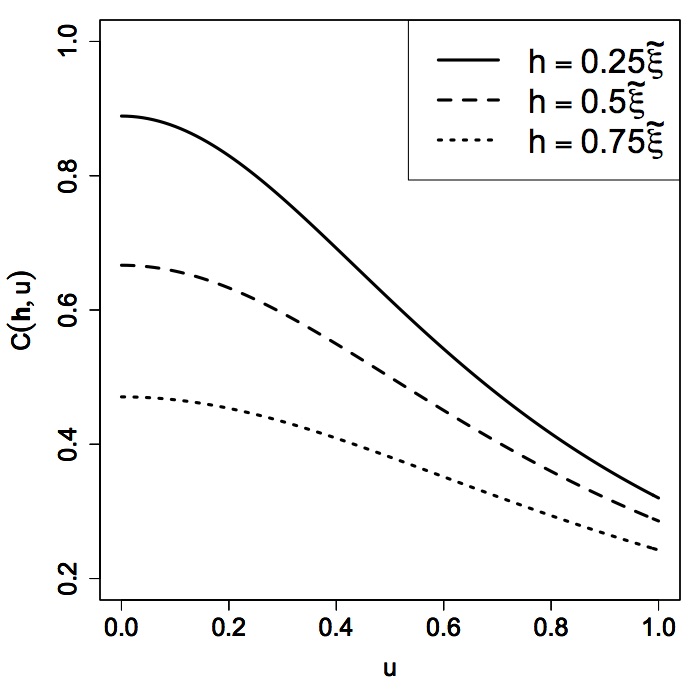}
\caption{The mapping $u \mapsto C(\bm{h},u)$ for the transport covariance $C(\bm{h},u)=\frac{1}{2} \{ C_S(\bm{h}-u\bm{\xi}) + C_S(\bm{h}+u\bm{\xi}) \}$, with $\bm{\xi}=(1,1)^\top$ and  Cauchy spatial covariance $C_S$. We fix spatial lags  of the form  $\bm{h}=\kappa \bm{\xi}$ (left), where dimple effect occurs for $|\kappa|>1/\sqrt{6}$, as well as spatial lags proportional to $\widetilde{\bm{\xi}}=(1,-1)^\top$, where dimple effect does not occur. }
\label{not_isotropic}
\end{figure}

Next, we focus on some results that directly emerge from Theorem \ref{teorema1}. If we assume that $\bm{V}$ is uniformly distributed over the spherical shell $\mathbb{S}^{d-1}$ of $\R^d$, and thus symmetrically distributed, its characteristic function $\varphi_{V}$ has expression 
 $\varphi_V(\bm{\eta}) = \Omega_{d}(\|\bm{\eta}\|)$ (see \citet{DP13}), where 
 \begin{equation}
  \label{aporte3} 
  \Omega_d(z) = \Gamma(d/2)(z/2)^{-(d-2)/2} J_{(d-2)/2}(z), \qquad z\in\mathbb{R}, 
  \end{equation}
with  $\Gamma$ being the Gamma function and $J_\nu$  the Bessel function of the first kind of degree $\nu$ \cite{Abramowitz-Stegun:1965}. Elementary properties of Bessel functions show that $\mathcal{H}_V = - (1/d) I_d$, where $I_d$ is the $d\times d$ identity matrix. Therefore, the mapping $F$ defined through Equation (\ref{aporte}) admits expression $F(\bm{h}) = -(1/d)\Delta C_S(\bm{h})$,  where $\Delta= \nabla^\top \nabla$ denotes the Laplacian operator.  We have deduced the following.
\begin{cor} 
\label{aporte2}
Suppose that $C_S$ is continuous and has second order derivatives on $\mathbb{R}^d\setminus \{\bm{0}\}$. If $\bm{V}$ is uniformly distributed on $\mathbb{S}^{d-1}$, then the   covariance  (\ref{transporte}) has a dimple  if, and only if, there exists two sets  $\Lambda_1$  and $\Lambda_2$, as in Definition \ref{def_dimple}, such that   $\Delta C_S (\bm{h})<0$  for  $\bm{h}\in \Lambda_1\setminus \{\bm{0}\}$, and  $\Delta C_S (\bm{h})>0$  for $\bm{h}\in \Lambda_2$.  
\end{cor}

Let us now turn to the popular case where the function $C_S$ that generates the Lagrangian covariance is radially symmetric, so that, there exists a continuous mapping $\phi:[0,\infty) \to \R$ such that $C_{S}(\bm{h})= \phi(\|\bm{h}\|)$. The class of such functions $\phi$ is uniquely identified with that of scale mixtures of the function $\Omega_d$ in Equation (\ref{aporte3}) and we refer the reader to \citet{DP13} with the references therein for a more recent discussion about this representation. For such a case, we have that $\Delta C_S(\bm{h}) =  \phi''(\|\bm{h}\|) + (d-1)\phi'(\|\bm{h}\|)/\|\bm{h}\|$, for all $\bm{h}\neq \bm{0}$, and Theorem \ref{teorema1}  reads as follows.  
\begin{cor}\label{corolario2}
Let $\bm{V}$ be uniformly distributed on $\mathbb{S}^{d-1}$ and let $C_S:\R^d \to \R$ be an isotropic covariance function, so that $C_S(\cdot) = \phi(\|\cdot\|)$, with $\phi$ being continuous and twice differentiable on $(0,\infty)$. Then,  the covariance  (\ref{transporte})   has a dimple  if, and only if, there exists a  constant $L \geq 0$  such that the function
\begin{equation}
 \label{last_condition}
r \mapsto  \phi''(r) + \frac{d-1}{r}\phi'(r), \qquad r>0,
\end{equation}  
is negative for $r\leq L$ and positive for  $r>L$.  
\end{cor}

In particular,  under the hypothesis of Corollary \ref{corolario2}, the  transport covariance $C$ is radial in space (see \citet{heidok22193}). The following example confirms that the  occurrence  of the dimple effect  could depend on the dimension  of the spatial domain.

\begin{figure}
\centering
\includegraphics[scale=0.25]{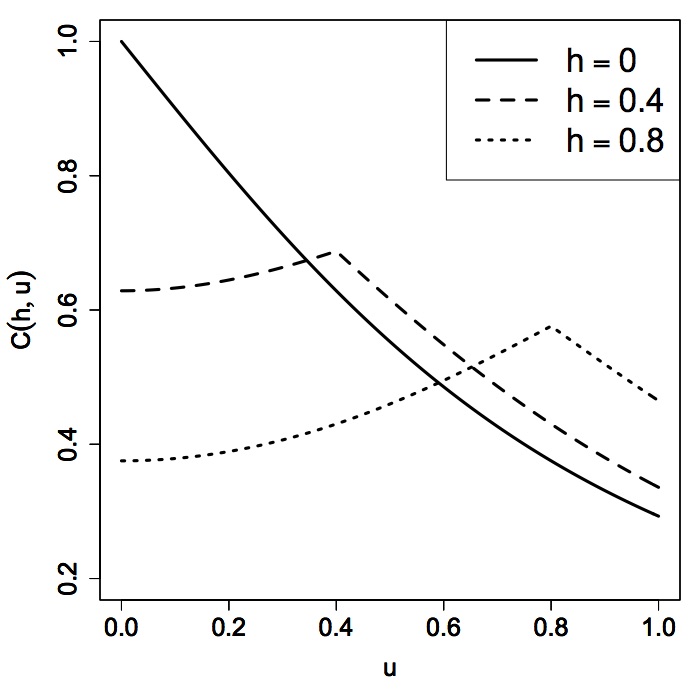} \hspace{2cm}  \includegraphics[scale=0.27]{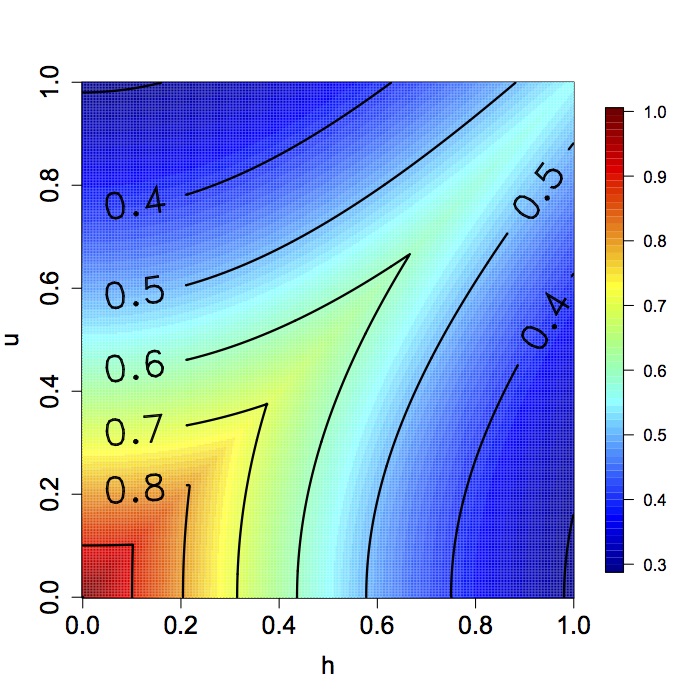}
\caption{The mapping $u \mapsto C(h,u)$ (left) for the transport covariance $C(h,u) = {\rm E}\{C_{S}(h- uV)\}$ on $\mathbb{R}\times \mathbb{R}$, with  $C_S$ in  the  Dagum family and $\Pr(V=\pm 1)=0.5$. We also provide the contour lines of $C$ on the unit square (right).}
\label{dagum}

\end{figure}

\begin{ejemplo}
Suppose that $\bm{V}$ is uniformly distributed on $\mathbb{S}^{d-1}$ and that $C_S$ belongs to the  Dagum family of  isotropic spatial  covariance functions \cite{berg2008,porcu2004}, that is,
\begin{equation*}
C_S(\bm{h}) = 1- \left(  \frac{\|\bm{h}\|^\gamma}{1+\|\bm{h}\|^\gamma} \right)^{\epsilon/\gamma},  \qquad \gamma\in (0,2], \quad \epsilon \in (0,\gamma).
\end{equation*}
Let  $\epsilon=1$ and $\gamma=2$. Then, mapping (\ref{last_condition}) is positively proportional to  $r\mapsto (4-d)r^2 + 1-d$.   Note that if $d\geq 4$, there is no dimple. Figure \ref{dagum} illustrates the dimple effect of the transport covariance for $d=1$. Since in this case the mapping (\ref{last_condition}) is strictly positive on $(0,\infty)$, the dimple arises immediately by taking any nonzero spatial lag.  
\end{ejemplo}

Additional  particular examples can be found in \citet{heidok22193}.  His illustrations   are consistent with our results.

\section{Transport effect models on spheres and dimple effect}
\label{sphere}

\subsection{General approach}
In this section, we consider space-time  random fields defined spatially over the  unit sphere   $\mathbb{S}^{d-1} \subset \mathbb{R}^d$. We give special emphasis to the cases $d=2,3$ representing the circle and the sphere of $\R^3$, respectively. This last case is especially important when modeling global data \cite{doi:10.1080/01621459.2015.1072541}.    For thorough studies on random fields on spheres, or spheres across time, we refer the reader to \cite{berg2016schoenberg,gneiting2013,marinucci2011random,doi:10.1080/01621459.2015.1072541}. In particular, the essay in \citet{gneiting2013} contains an impressive list of references, as well as an online supplement with a collection of open problems.

The formulation of  transport effect models over spheres across time should take into account the curvature of the sphere. We are not aware of any work related to such a construction over spheres, and proceed to illustrate a way to create a transport effect under this framework.   Let $\{Y(\bm{x}):\bm{x}\in\mathbb{S}^{d-1}\}$ be a Gaussian field on $\mathbb{S}^{d-1}$,  and let $\psi_S:[-1,1]\rightarrow \mathbb{R}$ be a  continuous function satisfying $\psi_S(\cos\theta) = \text{cov}\{Y(\bm{x}),Y(\bm{y})\}$, where $\theta: \mathbb{S}^{d-1} \times \mathbb{S}^{d-1}\to [0,\pi] $ denotes the great circle distance $\theta:= \theta(\bm{x},\bm{y}) = \arccos(\bm{x}^\top\bm{y})$, for $\bm{x},\bm{y}\in\mathbb{S}^{d-1}$. The class  of such functions is uniquely identified through Schoenberg representation (see \citet{gneiting2013} with the references therein).

  In order to construct a Lagrangian framework, we now suppose that the entire field $Y(\bm{x})$ moves in some random direction following the curvature of the sphere.  An appropriate way to represent   such a   displacement is through a  random rotation matrix $\cal{R}$ of order $d\times d$, that is, a random orthogonal matrix with determinant identically equal to 1. Recall that $\mathcal{R}$ is orthogonal if $\mathcal{R}\mathcal{R}^\top=\mathcal{R}^\top\mathcal{R}=I_d$, or simply $\mathcal{R}^{-1}=\mathcal{R}^\top$. The notion of randomness of $\mathcal{R}$ depends on the dimension of the sphere where the field is defined. For instance, if $d=2$, we could take two opposite directions of rotation, given by the clockwise and anti-clockwise movements.

All rotation matrices are diagonalizable  over the field of the complex numbers, namely $\mathcal{R} = QDQ^{-1}$, where the $i$-th column of $Q$ correspond to the $i$-th eigenvector of $\mathcal{R}$ and $D=\text{diag}\{\lambda_1,\hdots,\lambda_d\}$, with $\lambda_1,\hdots,\lambda_d$ being the eigenvalues of $\mathcal{R}$. Here, each eigenvalue admits expression $\lambda_i = \exp(\bm{\imath} \kappa_i)$, for some real constant $\kappa_i$, for $i=1,\hdots,d$. The last representation allows to define the powers of the matrix as follows  (see \citet{gantmacher1960theory}  for a complete discussion about  functions of matrices)
\begin{equation*}
\mathcal{R}^u :=  Q \begin{bmatrix}  \exp(\bm{\imath} \kappa_1 u) & \hdots & 0\\
                               \vdots  & \ddots & \vdots \\
                               0 & \hdots &  \exp(\bm{\imath} \kappa_d u)\\
  \end{bmatrix}	Q^{-1}, \qquad u\in\mathbb{R}.
\end{equation*}

We can finally define a   space-time field $Z(\bm{x},t)$ with Lagrangian dynamic on the sphere through the identity 
\begin{equation}
\label{lagrangian}
Z(\bm{x},t) =  Y(\mathcal{R}^t \bm{x}), \qquad (\bm{x},t)\in\mathbb{S}^{d-1}\times \mathbb{R}.
\end{equation}
The covariance function  corresponding to  (\ref{lagrangian}) is given by $K(\bm{x},\bm{y},t,t+u) = \text{cov}\{Z(\bm{x},t),Z(\bm{y},t+u)\}$, where
\begin{equation}\label{cov_sphere}
K(\bm{x},\bm{y},t,t+u) =  {\rm E}\{\psi_S( \bm{x}^\top \mathcal{R}^{u}  \bm{y} )\},  \qquad \bm{x},\bm{y} \in \mathbb{S}^{d-1}, t,u \in \R,
\end{equation}
and expectation is taken with respect to the random elements of $\mathcal{R}$.  

\begin{figure}
\centering{
\includegraphics[scale=0.2]{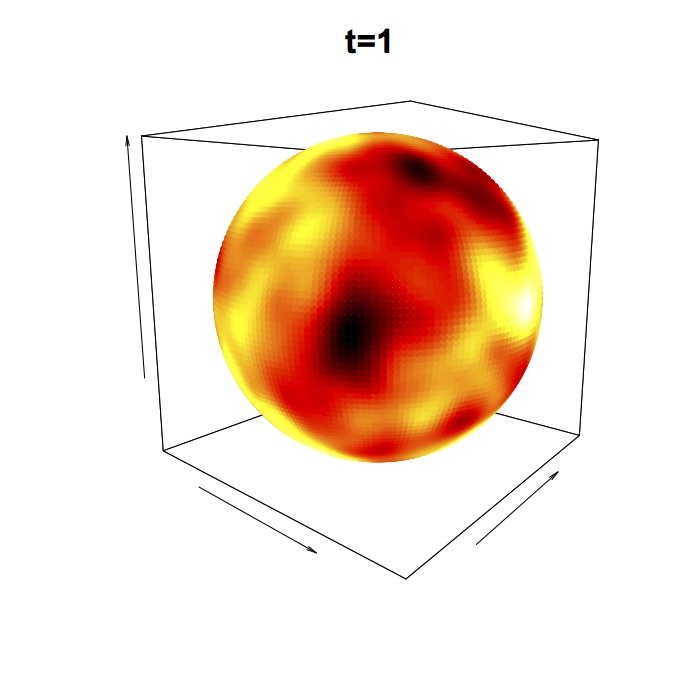} \includegraphics[scale=0.2]{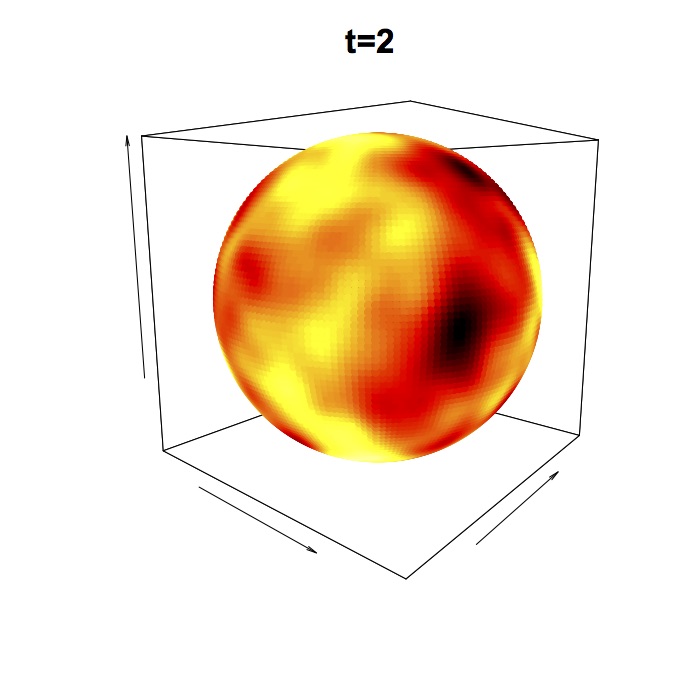} \includegraphics[scale=0.2]{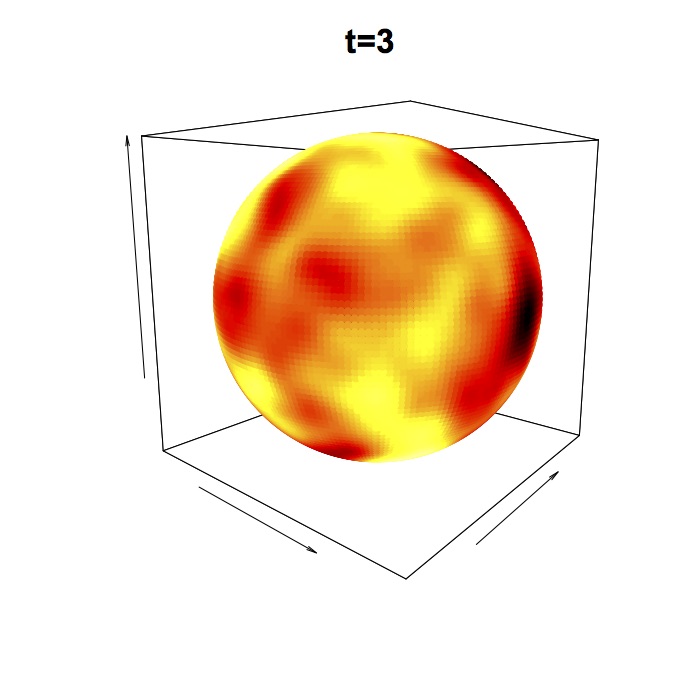}
}
\caption{Transport effect on $\mathbb{S}^2\times \R$, with $\psi_S$ in the multiquadric family. The field moves around the axis $(0,0,1)^\top$.}
\label{sim_s2}
\end{figure}

\begin{ejemplo}
Figure \ref{sim_s2}  illustrates the realization of a  field on $\mathbb{S}^2\times\R$,  generated from the multiquadric family 
 \begin{equation} 
 \label{multiquadric}
\psi_S(\cos\theta)   =   \frac{(1-\delta)^{2\tau}}{(1+\delta^2-2\delta \cos\theta)^\tau},  \qquad \delta\in(0,1), \tau>0,
\end{equation}
  which is a valid model on $\mathbb{S}^{d-1}$ for any dimension $d\geq 2$ \cite{gneiting2013}. The entire field  moves time-forward around the axis $(0,0,1)^\top$.  We set  $\delta=0.3$ and $\tau=0.5$, and consider three temporal instants $t=1,2,3$. 
\end{ejemplo}

In the following, we pay  attention to special cases of (\ref{cov_sphere}) that lead to spatial isotropy. More precisely, when there exists a function $\psi(\theta,u):  [0,\pi]\times \mathbb{R} \rightarrow \mathbb{R}$  such that $\psi(\theta,u) = K(\bm{x},\bm{y},t,t+u)$, for $\theta=\theta(\bm{x},\bm{y})$.  In addition, we focus on   mappings $\psi$ being symmetrical in the temporal argument.  We shall study the dimple effect on spheres under these conditions.  We  now  adapt  Definition \ref{def_dimple2} to the framework discussed in this section.

\begin{definition}	\label{def_dimple_sphere}
Consider a spatially isotropic and temporally symmetric covariance function $\psi(\theta,u)$, for $(\theta,u)\in [0,\pi]\times \mathbb{R}$,  associated to a random field on $\mathbb{S}^{d-1}\times \mathbb{R}$. We say that $\psi$ has a dimple in the temporal lag $u$ if there exists a  constant $L\in [0,\pi)$ such that
\begin{itemize}
\item[(1)] For fixed $\theta_0 \leq  L$,  the mapping $u \mapsto\psi(\theta_0,u)$ has a local maximum at $u=0$.
\item[(2)] For fixed $\theta_0 >  L$,  the mapping $u \mapsto \psi(\theta_0,u)$ has a local minimum at $u=0$.  
\end{itemize}
\end{definition}
Next, we study the cases $d=2$ and $d=3$.

\subsection{The simplest case: the unit circle}

First, we consider the unit circle  by   setting $d=2$ and we suppose that  two fields  move in opposite directions (clockwise and anti-clockwise). In this case, the rotation matrix has the form 
$$\mathcal{R}(\alpha) = \begin{bmatrix} \cos \alpha & -\sin\alpha\\ \sin \alpha & \cos \alpha  \end{bmatrix}, \qquad \alpha \in (0,2\pi].$$
Note that the relation $\mathcal{R}^u(\alpha) = \mathcal{R}(u\alpha)$ is satisfied for all $u\in\mathbb{R}$.   Fixing an angle $0<\alpha\leq 2\pi$, the matrices $\mathcal{R}(u\alpha)$ and $\mathcal{R}(-u\alpha)$,   describe the two opposite movements mentioned above by an angle proportional to $\alpha$.  If both directions have the same probability, the covariance function in Equation (\ref{cov_sphere})  reduces to
\begin{equation}\label{circle_cov}
K(\bm{x},\bm{y},t,t+u) =  \frac{1}{2}  \bigg\{ \psi_S( \bm{x}^\top \mathcal{R}(u \alpha) \bm{y} )  +  \psi_S( \bm{x}^\top \mathcal{R}(-u\alpha) \bm{y} )  \bigg\}, \qquad \bm{x},\bm{y}\in\mathbb{S}^1, t,u\in\mathbb{R}.
\end{equation}
Clearly, (\ref{circle_cov}) is symmetric in the temporal lag. Moreover, it depends on $\bm{x}$ and $\bm{y}$ only through their  great circle distance $\theta$.   In fact, direct inspection shows that $\bm{x}^\top \mathcal{R}(\alpha u) \bm{y} = \cos (\alpha u) \bm{x}^\top  \bm{y}   - \sin (\alpha u) (x_1y_2-x_2y_1)$. In addition,  $\bm{x}^\top  \bm{y}  = \cos\theta$ and $  (x_1y_2-x_2y_1)^2 = 1 -  (\bm{x}^\top  \bm{y})^2  = \sin^2\theta$.  Thus,   $\bm{x}^\top \mathcal{R}(\alpha u) \bm{y}$ coincides with   $\cos(\theta + \alpha u)$  or with  $\cos(\theta - \alpha u)$. In conclusion, we can derive the following representation for (\ref{circle_cov})
\begin{equation}  \label{circle_cov2}
\psi(\theta,u) =  \frac{1}{2}\bigg \{ \psi_S( \cos(\theta + u\alpha) )  +  \psi_S( \cos(\theta - u\alpha))  \bigg\}, \qquad \theta \in [0,\pi], u\in\mathbb{R}.
\end{equation}
This expression is the crux for the characterization of the dimple on the unit circle.

\begin{theorem} \label{teorema_circle}
Let $\psi_S$ be continuous and  twice differentiable on $(-1,1)$. The transport covariance (\ref{circle_cov2}) has a dimple along the temporal lag  if, and only if, there exists a constant $L \in[0,\pi)$ such that the function
\begin{equation} 
\label{condition_circle}
\theta  \mapsto   \psi_S^{''}(\cos\theta) \sin^2\theta   - \psi_S^{'}(\cos\theta) \cos\theta, \qquad \theta \in (0,\pi),
\end{equation}  
is negative for $\theta \leq L$ and positive for $\theta>L$.
\end{theorem}

We do not report the proof of Theorem \ref{teorema_circle}, since an easy calculation shows that at $u=0$ the partial derivative $\partial \psi(\theta,u)/\partial u$ is zero and   $\partial^2 \psi(\theta,u)/\partial u^2$ has the same sign as (\ref{condition_circle}).

\begin{figure}
\begin{center}
\includegraphics[scale=0.25]{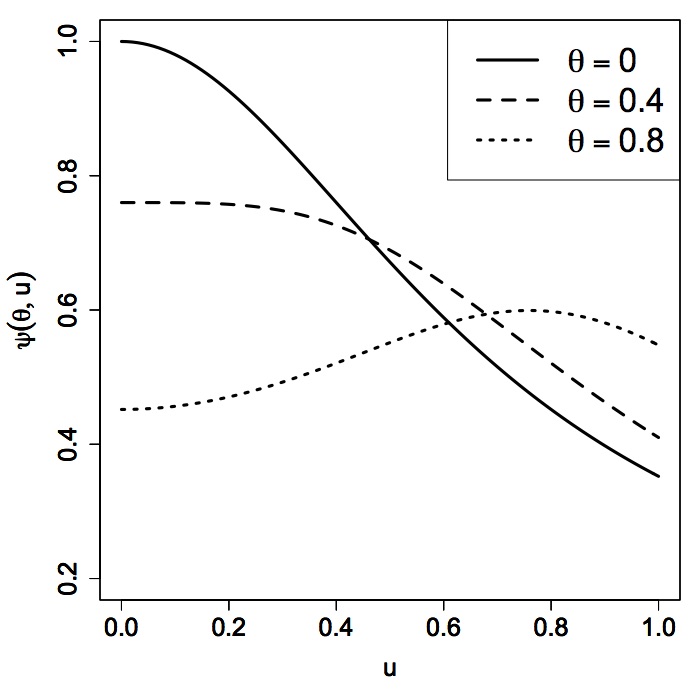}  \hspace{2cm} \includegraphics[scale=0.27]{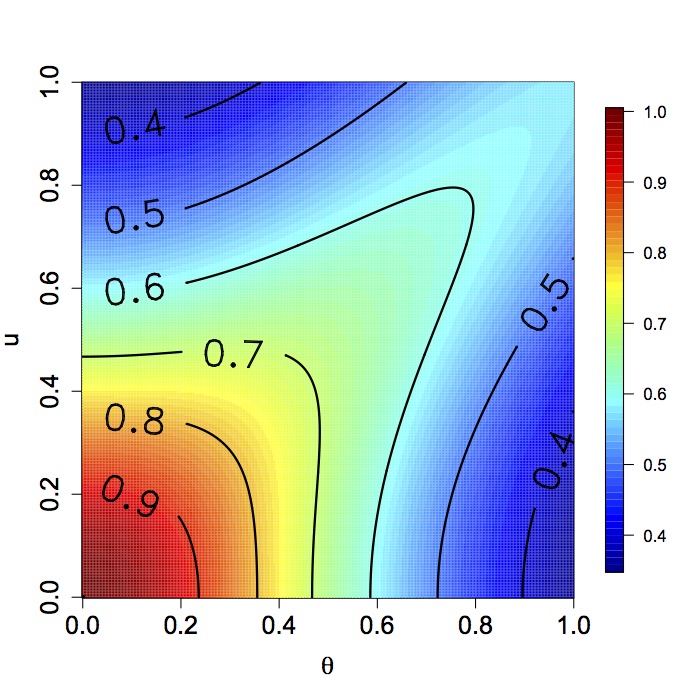}
\caption{The mapping $u \mapsto\psi(\theta,u)$ (left) for the transport covariance $\psi(\theta,u) =\{ \psi_S( \cos(\theta + u\alpha) )  +  \psi_S( \cos(\theta - u\alpha)) \}/2$ on $\mathbb{S}^1\times \mathbb{R}$, with  $\psi_S$ in  the  multiquadric family. We also provide the contour lines of $\psi$ on the unit square (right).}
\label{mult}
\end{center}
\end{figure}

\begin{ejemplo}
Consider the multiquadric model (\ref{multiquadric}) on $\mathbb{S}^1$, with $\tau = 1$ and $\delta \in (0,1)$. Then, (\ref{condition_circle}) has the same sign as   $2\delta\cos^2\theta  +  (1+\delta^2) \cos\theta  - 4\delta$.  Therefore, the problem is reduced to obtain the roots of this quadratic equation in $\cos\theta$. In fact, denote as $\zeta_1=(-(1+\delta^2) + \sqrt{(1+\delta^2)^2+32\delta^2}) / (4\delta)$ and $\zeta_2=(-(1+\delta^2) - \sqrt{(1+\delta^2)^2+32\delta^2}) / (4\delta)$ such roots, which satisfy that $0 < \zeta_1 <1$ and $\zeta_2 < -2$, for all $0<\delta<1$. Therefore, we have a dimple by taking $L=\arccos\left( \zeta_1\right)$ in Theorem \ref{teorema_circle}. Figure \ref{mult} depicts the dimple effect considering $\alpha=1$ and $\delta=0.5$.
\end{ejemplo}

\subsection{Transport effects over the Earth's surface: the case  $\mathbb{S}^2\times \R$}

The \textit{Rodrigues rotation formula} \cite{Kuipers:710564} establishes that  a rotation in $\mathbb{R}^3$ by an angle $\alpha$, with respect to an arbitrary axis determined by the unit vector $\bm{\omega}=(\omega_1,\omega_2,\omega_3)^\top\in\mathbb{R}^3$,  can be written as
\begin{equation*}
\mathcal{R}_{\bm{\omega}}(\alpha) = W \sin(\alpha)  +  (I_3 - \bm{\omega}\bm{\omega}^\top) \cos(\alpha) + \bm{\omega}\bm{\omega}^\top,
\end{equation*} 
where 
\begin{equation*} 
W = \begin{bmatrix}
0 & -\omega_3 & \omega_2\\
\omega_3 & 0 & -\omega_1\\
-\omega_2 & \omega_1 & 0
\end{bmatrix}.
\end{equation*}
The rotation matrix satisfies the relation $\mathcal{R}^u_{\bm{\omega}}(\alpha) = \mathcal{R}_{\bm{\omega}}(u \alpha)$, for all $u\in\mathbb{R}$. Therefore, if we fix an angle $0 < \alpha \leq 2\pi$ and  take  the axis $\bm{\omega}$ randomly,  the transport model is given by 
\begin{equation} 
\label{transport_s2}
K(\bm{x},\bm{y},t,t+u) = {\rm E}\{ \psi_S(   \bm{x}^\top  \mathcal{R}_{\bm{\omega}}(u \alpha) \bm{y}   ) \},  \qquad  \bm{x},\bm{y}\in\mathbb{S}^2, t,u\in\mathbb{R}.
\end{equation}

The following result establishes a sufficient condition on the distribution of $\bm{\omega}$ such that  (\ref{transport_s2}) is spatially isotropic.

\begin{proposition}
\label{prop_isotropy}
Let $\bm{\omega}$ be uniformly distributed on $\mathbb{S}^2$. Then, there exists a function $\psi:[0,\pi]\times\R \to \R$ such that the covariance in Equation (\ref{transport_s2}) can be written as $K(\bm{x},\bm{y},t,t+u)= \psi(\theta,u)$, where $\theta$ is the great circle distance between $\bm{x}$ and $\bm{y}$ in $\mathbb{S}^2$.
\end{proposition}

The proof of Proposition \ref{prop_isotropy} is deferred to Appendix \ref{app2}.  In addition, it is easy to show that, for $\bm{\omega}$  uniformly distributed on $\mathbb{S}^2$, the covariance $\psi(\theta,u)$  is symmetric in the temporal argument. In fact, note that
$$  \psi(\theta, -u) = {\rm E}\{ \psi_S(   \bm{x}^\top  \mathcal{R}_{\bm{\omega}}(-u \alpha) \bm{y}   ) \} = {\rm E}\{ \psi_S(   \bm{x}^\top  \mathcal{R}_{-\bm{\omega}}(u \alpha) \bm{y}   ) \} = \psi(\theta,u),$$
where the last equality comes from the fact that $-\bm{\omega}$ is uniformly distributed on $\mathbb{S}^{2}$.

\begin{theorem}
\label{teo_sphere}
Let $\psi_S$ be continuous and  twice differentiable on $(-1,1)$ and  the axis of rotation $\bm{\omega}$ being uniformly distributed on $\mathbb{S}^2$.  Then, the covariance $\psi(\theta,u)={\rm E}\{ \psi_S(   \bm{x}^\top  \mathcal{R}_{\bm{\omega}}(u \alpha) \bm{y}   ) \}$ has a dimple in the temporal lag  if, and only if, there exists a constants $L\in [0,\pi)$ such that the function
\begin{equation} 
\label{condition_s2}
\theta  \mapsto   \psi_S^{''}(\cos\theta) \sin^2\theta   -  2 \psi_S^{'}(\cos\theta) \cos\theta, \qquad \theta \in (0,\pi),
\end{equation}  
is negative for $\theta\leq L$ and positive for $\theta>L$.
\end{theorem}

We refer the reader to Appendix \ref{app3} for a proof of Theorem \ref{teo_sphere}.

\begin{figure}
\begin{center}
\includegraphics[scale=0.25]{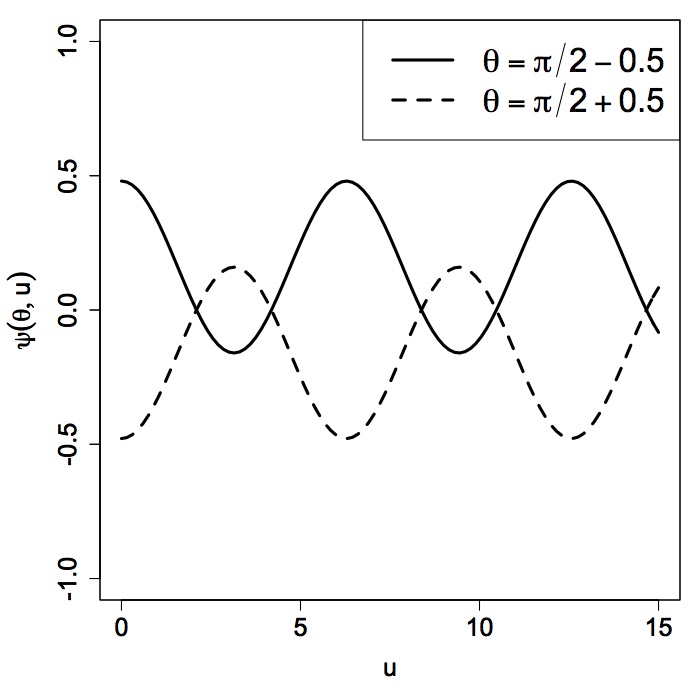} \hspace{2cm} \includegraphics[scale=0.27]{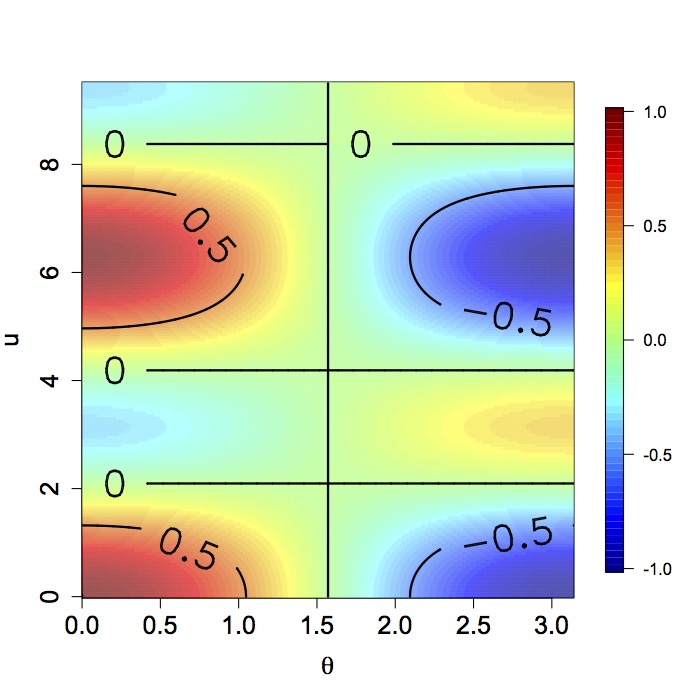}
\caption{The mapping $u \mapsto\psi(\theta,u)$ (left) for the transport covariance $\psi(\theta,u)={\rm E}\{ \psi_S(   \bm{x}^\top  \mathcal{R}_{\bm{\omega}}(u \alpha) \bm{y}   ) \} $ on $\mathbb{S}^2\times \mathbb{R}$, with  $\psi_S=\cos\theta$, and two fixed  values for $\theta$. We also provide the contour lines of $\psi$ on $[0,\pi]\times [0,9.5]$ (right).}
\label{dimple_esfera}
\end{center}
\end{figure}

\begin{ejemplo}
We start with the covariance $\psi_S(\cos\theta) = \cos\theta$ on  $\mathbb{S}^2$. The associated transport model is $\psi(\theta,u) = \{1+2\cos(u\alpha)\} (\cos\theta)/3$. The  mapping (\ref{condition_s2}) is given by $\theta \mapsto -2\cos\theta$, which is negative on $(0,\pi/2)$ and positive on $(\pi/2,\pi)$.  Figure \ref{dimple_esfera} illustrates the dimple effect for this model with $\alpha=1$.
\end{ejemplo}

\section{Concluding remarks}

This paper studied the occurrence of dimple in space-time  geostatistical models coming from transport effects. In our formulation, the spatial domain  can be  either $\mathbb{R}^d$ or $\mathbb{S}^d$. We provided characterization theorems and several examples to illustrate our findings.  We showed that dimple effects can  occur in situations where there are no prevailing  directions, and that its presence is strongly connected to the spatial dimension of the problem. 

We believe that the Lagrangian framework as well as the definition of dimple can be generalized in order to cover  non stationary behaviors in both space and time.     An interesting approach is to consider wind directions  having a dynamically updated distribution according to the state of the atmosphere.  In this case, the covariance would not be stationary in time and our results could  be extended to this more general scenario.   On the other hand, axially symmetric covariances are an interesting, and physically motivated,  alternative to model environmental processes on the globe, allowing for  non stationarity in the latitude component. The formulation of transport effects, and the study of the corresponding dimple effect, cannot be directly extended to this context, since arbitrary rotations  could  alter the underlying spatial structure. Indeed, more effort is needed to cover this framework.

 It might be worth noticing that we have proposed a simple model for transport effect in order to explore the presence of a dimple. More general and sophisticated models have been proposed and we refer to the book by \citet[pp. 319-320]{cressie2011statistics}. For future researches, a big challenge would be to explore the dimple problem for, e.g., non homogeneous transport effects. 

Another important aspect is that the models proposed in this paper are radially symmetric in Euclidean spaces and isotropic over spheres. Such an assumption is of course limited for the analysis of real data. Quoting Noel Cressie (personal communication):

\textit{I must say that, in my experience, environmental phenomena almost never follow this paradigm,
even when mean processes are removed by working with anomalies. Global geophysical phenomena
are generally highly non-stationary in both space and time.} 

At the same time, isotropic models are important because they can be used as building blocks for more sophisticated models, as suggested by \citet{porcu2017modeling}.

\section*{Acknowledgments}
The authors are very grateful  to Mandy Hering and Noel Cressie   for fruitful discussions and suggestions during the preparation of the manuscript.   Alfredo Alegr\'ia is supported by \textit{Beca CONICYT-PCHA/Doctorado Nacional/2016-21160371} and  \textit{Programa de Incentivos a la Iniciaci\'on Cient\'ifica  2016, University Federico Santa Mar\'ia}. Emilio Porcu is supported by \textit{Proyecto FONDECYT Regular No 1130647}.

\appendix

\renewcommand*{\thesection}{\Alph{section}}

\section*{Appendices}

\section{Proof of Theorem \ref{teorema1}}
\label{app1}

Before we state the proof of Theorem  \ref{teorema1}, we need some preliminary results.  The main ingredients needed for the results following subsequently rely on Bochner's characterization \cite{bochner1955harmonic} of continuous covariance functions on $\R^d$ as being the Fourier transforms of positive and bounded measures $\mu: \R^d \to \R$, also called spectral measures:
\begin{equation}
\label{bochner}
C_S(\bm{h}) = \int_{\mathbb{R}^d} \exp ( \bm{\imath} \bm{h}^\top\bm{\omega} ) \mu (\text{d}\bm{\omega}), \qquad \bm{h} \in \R^d.
\end{equation}
Coupling Bochner's representation with the construction in Equation (\ref{transporte}), we obtain a useful lemma being needed for the proof of the main result.
\begin{lemma} \label{lemmata}
Let $C: \R^d \times \R \rightarrow \R$ be the function defined through Equation (\ref{transporte}). Then, the following identity is true:
\begin{equation}\label{spec_repres}
C(\bm{h},u) = \int_{\mathbb{R}^d}  \exp ( \bm{\imath}  \bm{h}^\top\bm{\omega}  ) \varphi_V(-u \bm{\omega} ) \mu(\text{d}\bm{\omega}), \qquad (\bm{h},u) \in \R^d \times \R,
\end{equation}
where  $\mu$ is the spectral measure of $C_S$, defined according to Equation (\ref{bochner}).
\end{lemma}

\textbf{Proof of Lemma \ref{lemmata}.}
Using the spectral representation of $C_S$ and Fubini's theorem we have 
\begin{eqnarray*}
C(\bm{h},u)      =      {\rm E}\left(   \int_{\mathbb{R}^d} \exp\{ \bm{\imath} (\bm{h}-u\bm{V})^\top\bm{\omega}\}  \mu(\text{d}\bm{\omega}) \right)  =  \int_{\mathbb{R}^d} \exp (\bm{\imath} \bm{h}^\top\bm{\omega} )    {\rm E} \{  \exp ( -\bm{\imath}u \bm{V}^\top  \bm{\omega} )  \} \mu(\text{d}\bm{\omega}).     
\end{eqnarray*}  Finally, note that $\varphi_V(-u \bm{\omega} )= {\rm E} \{  \exp(-\bm{\imath} u \bm{V}^\top \bm{\omega})  \} $. These facts complete the proof.

\textbf{Proof of Theorem \ref{teorema1}.}
 We give a constructive proof. Consider the spatial lag $\bm{h}\neq\bm{0}$. The representation  (\ref{spec_repres}) in the previous Lemma and direct inspection show that $\partial C(\bm{h},u)/\partial u = - \int_{\mathbb{R}^d} \exp( \bm{\imath}  \bm{h}^\top \bm{\omega}  )      \{ \bm{\omega}^\top  \nabla \varphi_V(-u\bm{\omega}) \} \mu(\text{d}\bm{\omega})$, where the exchange of derivative with the integral is justified by the fact that the integrand is uniformly bounded.
Since $\bm{V}$ is symmetrically distributed, we deduce that $ \nabla \varphi_V(\bm{0})=\bm{0}$, implying 
 $\partial C(\bm{h},u)/\partial u  |_{u=0} = 0 $. Moreover, direct inspection shows that
\begin{equation*}
\frac{\partial^2 C(\bm{h},u)}{\partial u^2}   \Big|_{u=0}   =  \int_{\mathbb{R}^d} \exp( \bm{\imath} \bm{h}^\top\bm{\omega} )     (\bm{\omega}^\top  \mathcal{H}_V \bm{\omega})    \mu(\text{d}\bm{\omega}). 
\end{equation*}
On the other hand,   the Hessian matrix of $C_S$ is given by
\begin{equation*}
\nabla^2  C_S(\bm{h})  =  - \int_{\mathbb{R}^d} \exp ( \bm{\imath}  \bm{h}^\top \bm{\omega} )   (\bm{\omega} \bm{\omega}^\top) \mu  (\text{d}\bm{\omega}),  \qquad \bm{h} \in \R^d,
\end{equation*}
where integration is taken componentwise.
Therefore, the following equality holds 
\begin{eqnarray*}
   F(\bm{h})  &     =     &  -   \int_{\mathbb{R}^d} \exp ( \bm{\imath}  \bm{h}^\top\bm{\omega} )  \tr( \mathcal{H}_V  \bm{\omega} \bm{\omega}^\top  ) \mu  (\text{d} \bm{\omega}),
\\
  &    =     &  -  \int_{\mathbb{R}^d} \exp ( \bm{\imath}  \bm{h}^\top\bm{\omega} )  \tr( \bm{\omega}^\top  \mathcal{H}_V  \bm{\omega}  )  \mu(\text{d}\bm{\omega}), \\
 &  =  &   - \frac{\partial^2 C(\bm{h},u)}{\partial u^2}   \Big|_{u=0}. 
   \end{eqnarray*} 
Finally, for any fixed $\bm{h}\neq \bm{0}$, the function $u \mapsto C(\bm{h},u)$ has a  local minimum or maximum at the origin depending on the sign of $F(\bm{h})$. The proof is completed.

\section{Proof of Proposition \ref{prop_isotropy}}
\label{app2}

Let $\bm{x}$ and $\bm{y}$ be linearly independent points on $\mathbb{S}^{2}$. We need to show that if $\bm{\omega}$ is uniformly distributed on $\mathbb{S}^2$, the distribution of the quadratic form 
\begin{eqnarray*}
 \bm{x}^\top  \mathcal{R}_{\bm{\omega}}(u \alpha) \bm{y}     &    =      &     \cos(u\alpha)  (\bm{x}\top\bm{y})  +\sin(u\alpha)   (\bm{x}^\top W \bm{y})    + \{1 -\cos(u\alpha) \} (\bm{x}^\top \bm{\omega} ) (\bm{\omega}^\top \bm{y} ) \\
 &  =  &    \cos(u\alpha) \cos\theta - \sin(u\alpha)     \{\bm{\omega}^\top (\bm{x}\times\bm{y}) \}          + \{1 -\cos(u\alpha)\}  (\bm{x}^\top \bm{\omega})  ( \bm{y}^\top \bm{\omega} ),
\end{eqnarray*}
depends on $\bm{x}$ and $\bm{y}$ only through their great circle distance $\theta$, where $\times$ denotes the cross product in $\mathbb{R}^3$.   In fact, the rotation invariance property of $\bm{\omega}$ allows to consider the spherical coordinates with respect to an arbitrary orthonormal basis in $\mathbb{R}^3$. We take the mutually orthogonal  axis $\bm{e}_1' = \bm{x}$,  $\bm{e}_2'= ( \bm{x}\times\bm{y}) / \sin\theta $, and  $\bm{e}_3'= \{ \bm{x}\times(\bm{x}\times\bm{y}) \}/\sin\theta$. Moreover, the properties of cross product imply that $\bm{e}_3' =\{(\cos\theta) \bm{x} - \bm{y}\}/\sin\theta $. Then, $\bm{\omega}$ can be written with respect to this basis as
$$ \bm{ \omega } =    (\bm{\omega}^\top \bm{x})    \bm{e}_1'  +    \frac{ 1}{\sin\theta} \{\bm{\omega}^\top (\bm{x}\times\bm{y}) \}  \bm{e}_2'    +  \frac{ 1}{\sin\theta}  \{ (\cos\theta)(\bm{\omega}^\top  \bm{x}) - (\bm{\omega}^\top \bm{y}) \}\bm{e}_3'.$$
Let $0\leq \phi_1 \leq \pi$ and $0 \leq  \phi_2 < 2\pi$ be the azimuth and polar angles, respectively. Then, we have
\begin{eqnarray*}
   \bm{\omega}^\top \bm{x}    & = &  \cos \phi_1, \\
\frac{ \bm{\omega}^\top (\bm{x}\times\bm{y}) }{\sin\theta}     & = & \sin\phi_1 \sin\phi_2, \\
 \frac{   (\cos\theta)(\bm{\omega}^\top  \bm{x}) - (\bm{\omega}^\top \bm{y})  }{\sin\theta}  & = & \sin\phi_1 \cos\phi_2.
\end{eqnarray*}
Therefore, the transport covariance on  spheres across time can be represented through the following integral
\begin{multline*}
\psi(\theta,u) = \frac{1}{4\pi} \int_0^{2\pi}\int_0^\pi   \psi_S \bigg(        \cos(u\alpha) \cos\theta - \sin(u\alpha)   \sin\phi_1 \sin\phi_2 \sin\theta   
\\ + \cos\phi_1 \{1 -\cos(u\alpha)\}  \{  \cos\phi_1 \cos\theta - \sin\phi_1\cos\phi_2 \sin\theta   \}   \bigg)   \sin\phi_1 \text{d}\phi_1 \text{d}\phi_2.
\end{multline*} The proof is completed.

\section{Proof of Theorem \ref{teo_sphere}}
\label{app3}

 We have the quadratic form $\bm{x}^\top  \mathcal{R}_{\bm{\omega}}(u \alpha) \bm{y}       =        (\bm{x}^\top  W \bm{y})   \sin(u\alpha)        +          \{\bm{x}^\top  (I_3 - \bm{\omega}\bm{\omega}^\top) \bm{y}  \}  \cos(u\alpha)       +           (\bm{x}^\top  \bm{\omega})  (\bm{\omega}^\top \bm{y}).$ Thus, for $\theta\neq 0$, differentiation   respect to $u$ gives
\begin{eqnarray}
\frac{\partial \psi(\theta,u) }{\partial u}   & = &    {\rm E}\left( \psi^{'}_S(   \bm{x}^\top  \mathcal{R}_{\bm{\omega}}(u \alpha) \bm{y}   )  \frac{\partial \{ \bm{x}^\top  \mathcal{R}_{\bm{\omega}}(u \alpha) \bm{y}\}}{\partial u}   \right) \nonumber \\
   &     =    &     \alpha  {\rm E}\left( \psi^{'}_S(   \bm{x}^\top  \mathcal{R}_{\bm{\omega}}(u \alpha) \bm{y}   )  \bigg[     (\bm{x}^\top W \bm{y} )  \cos(u\alpha)        -         \{ \bm{x}^\top  (I_3 - \bm{\omega}\bm{\omega}^\top) \bm{y} \}   \sin(u\alpha)    \bigg]  \right),  \nonumber
\end{eqnarray}
so that 
$$\frac{\partial\psi(\theta,u)}{  \partial u } \big|_{u=0}    =    \alpha   \psi^{'}_S(\cos \theta)  {\rm E}\left(     \bm{x}^\top W \bm{y}    \right)   =  - \alpha   \psi^{'}_S(\cos \theta)     {\rm E}\left(             \bm{\omega}^\top \frac{ \bm{x}\times\bm{y} }{ \|\bm{x}\times\bm{y}\| }           \right)     \sin\theta  =0.$$  
On the other hand,
\begin{multline*}
\frac{\partial^2 \psi(\theta,u) }{\partial u^2}   =   \alpha^2\bigg(    
 {\rm E}\bigg\{ \psi^{''}_S(   \bm{x}^\top  \mathcal{R}_{\bm{\omega}}(u \alpha) \bm{y}   )  \bigg[     (\bm{x}^\top  W \bm{y})   \cos(u\alpha)        -        \{\bm{x}^\top  (I_3 - \bm{\omega}\bm{\omega}^\top) \bm{y}\}    \sin(u\alpha)    \bigg]^2       \\
  +     \psi^{'}_S(   \bm{x}^\top  \mathcal{R}_{\bm{\omega}}(u \alpha) \bm{y}   )  \bigg[   -  (\bm{x}^\top W \bm{y})   \sin(u\alpha)        -       \{  \bm{x}^\top  (I_3 - \bm{\omega}\bm{\omega}^\top) \bm{y} \}   \cos(u\alpha)     \bigg]     \bigg\}
\bigg).
\end{multline*}
Then,
\begin{eqnarray*}
\frac{\partial^2\psi(\theta,u) }{\partial u^2}  \Big|_{u=0}     &   =   &         \alpha^2\bigg(      
  \psi^{''}_S( \cos \theta  )       {\rm E}\{    (\bm{x}^\top  W \bm{y})^2   \} 
  -     \psi^{'}_S(  \cos \theta  )        [ \bm{x}^\top  \{I_3 -   {\rm E}(\bm{\omega}\bm{\omega}^\top)\}\bm{y}   ]
\bigg)    \\
&   =   &         \alpha^2\bigg(      
  \psi^{''}_S( \cos \theta  )       {\rm E}\left\{     \left(\bm{\omega}^\top  \frac{ \bm{x}\times\bm{y} }{ \|\bm{x}\times\bm{y}\| }  \right)^2   \right\} 
\sin^2 \theta  -   \frac{2}{3}  \psi^{'}_S(  \cos \theta  )       \cos\theta 
\bigg)       \\
&   =   &         \alpha^2\bigg(      
 \frac{1}{3} \psi^{''}_S( \cos \theta  )    \sin^2 \theta   
  -   \frac{2}{3}  \psi^{'}_S(  \cos \theta  )       \cos\theta 
\bigg) ,   
\end{eqnarray*}
where we have used  that ${\rm E}(\bm{\omega}\bm{\omega}^\top) = (1/3)I_3$ and  ${\rm E}\{(\bm{\omega}^\top\bm{e})^2\} = 1/3$ for any unit vector $\bm{e} \in\mathbb{R}^3$.

\section*{\refname}
\bibliographystyle{elsarticle-harv} 
\bibliography{mybib}

\end{document}